\documentclass{amsart}

\usepackage[english]{babel}
\usepackage{amssymb,amsmath,amsthm}

\newtheorem{theorem}{Theorem}
\newtheorem{proposition}{Proposition}
\newtheorem{lemma}{Lemma}
\theoremstyle{definition}
\newtheorem{definition}{Definition}

\newcommand{\modulus}[1]{\left\lvert #1 \right\rvert}
\newcommand{\norm}[1]{\left\| #1 \right\|}
\newcommand{\mklm}[1]{\left\{ #1 \right\}}
\newcommand{\eklm}[1]{\left\langle #1 \right\rangle}

\newcommand{\C}{{\mathbb C}}
\newcommand{\R}{{\mathbb R}}

\newcommand{\D}{{\mathcal D}}

\renewcommand{\H}{{\mathcal H}}

\newcommand{\M}{{\mathcal M}}
\newcommand{\T}{{\mathcal T}}
\renewcommand{\P}{{\rm P}}
\renewcommand{\O}{{\mathcal O}}

\newcommand{\Lint}{{\rm L}}
\newcommand{\Schw}{{\mathcal S}}

\newcommand{\hypP}{{\Gamma^+_P}}
\newcommand{\chypP}{{\Gamma^{(c)}_P}}

\renewcommand{\O}{{\mathrm O}}
\newcommand{\SU}{{\mathrm SU}}

\newcommand{\Linh}{{\mathrm P_+^\uparrow(3,1)}}

\newcommand{\Linhgen}{{\mathrm P(3,1)}}
\newcommand{\LC}{{\mathrm L_+(C)}}

\newcommand{\ad}{\text{ad}\,}

\renewcommand{\Im}{\mathrm{Im}\,}
\renewcommand{\Re}{\mathrm{Re}\,}
\DeclareMathOperator{\supp}{supp}

\DeclareMathOperator{\convh}{convh}

\begin{document}

\author{Pablo Ramacher}
\title[Localization of elementary systems in the theory of Wigner]{Modular localization of elementary systems in the theory of Wigner}
\address{Pablo Ramacher, Humboldt--Universit\"at zu Berlin, Institut f\"ur Reine Mathematik, Ziegelstr. 13a, D--10099 Berlin, Germany}
\email{ramacher@mathematik.hu-berlin.de}
\thanks{Supported by the SFB 288 of the DFG}

\begin{abstract}
Starting from Wigner's theory of elementary systems and following a recent approach of Schroer we define certain subspaces of localized wave functions in the underlying Hilbert space with the help of the theory of modular von--Neumann algebras of Tomita and Takesaki. We characterize the elements of these subspaces as boundary values of holomorphic functions in the sense of distribution theory and show that the corresponding holomorphic functions satisfy the sufficient conditions of the theorems of Paley--Wiener--Schwartz and H\"{o}rmander.
\end{abstract}

\maketitle

\section{Introduction}

The subject of this paper is the localization of elementary systems in the sense of Wigner. These are quantum mechanical systems whose states are all obtainable from any state by relativistic transforms and superposition. They constitute a relativistic invariant linear manifold and the corresponding wave functions satisfy relativistic invariant wave equations. 

Bargmann and Wigner \cite{bargmann-wigner} realized that these wave equations can be replaced by representations of the inhomogeneous Lorentz group given by these same  equations. A classification of all possible representations amounts then to a classification of all possible relativistic wave equations. In this way it is possible not only to construct solutions of the wave equations but also to specify their relevant invariant properties. It is natural to realize these representations in momentum space since the momenta and energies of the system, but not the coordinates, are defined by the Lorentz group as infinitesimal translations.

A priori it is not clear which localization properties do correspond to the different elementary systems, since the coordinates which appear as arguments in the coordinate space wave functions are not eigenvalues of the position operator conjugate to the momentum operator.

General quantum field theory  in the sense of Haag, Araki and Kastler \cite{haag} is primarily concerned with local operations. Thus to each open region in Minkowski space there is associated an algebra of operators acting on the underlying Hilbert space  which are interpreted as physical operations or observables that can be performed within this region. The states of the system are then defined as positive linear functionals over these algebras. Since it is sufficient to consider only bounded operators one is led to the study of von--Neumann algebras; their  properties can be analysed independently of the generating fields. By the theory of modular von--Neumann algebras of Tomita and Takesaki  \cite{tomita}, \cite{takesaki} it is possible to associate operators to certain states and space time regions which contain  important features of the theory. It was first pointed out by Schroer \cite{schroer1,schroer3} that by knowing these modular operators for certain regions one can associate real subspaces of localized wave functions in the original Wigner representation space to space time regions in Minkowski space. These real subspaces can then be used to construct free quantum theories. Recent work in this direction has also been done by Brunetti, Guido and Longo \cite{brunetti-guido-longo2,brunetti-guido-longo1}.

In this paper we characterize the elements of these subspaces as boundary values of analytic functions in the sense of distribution theory which fullfill certain boundary conditions and we show that these analytic functions are the Fourier--Laplace transforms of distributions with support in the considered closed, but not necessarily compact, convex  regions. We restrict ourselves to the case of the massive scalar field, but all our considerations can be carried over to arbitrary quantum fields.

\section{Representations of the inhomogeneous Lorentz group}

Let $\R^4$ be the four dimensional Minkowski space with coordinates $x^0,\, x^1,\,x^2,\,x^3$ and metric tensor $g$ given by $g_{00}=1, \,g_{11}=g_{22}=g_{33}=-1,\,g_{ij}=0,\, i\not=j$. The group of all linear transformations which leave the quadratic form
\begin{displaymath}
(x^0)^2-(x^1)^2-(x^2)^2-(x^3)^2
\end{displaymath}
invariant is called the general homogeneous Lorentz group. An inhomogeneous Lor\-entz transform is a transformation which consists of a homogeneous Lorentz transform together with a translation in Minkowski space, the translation being performed after the homogeneous Lorentz transform. The component of the unity of the general inhomogeneous Lorentz group is denoted by $\Linh$, the proper orthochronous inhomogeneous Lorentz group.

According to Wigner \cite{wigner1} the unitary irreducible representations of the inhomogeneous Lorentz group can be classified as follows:
\begin{theorem}
The representations of class $P^+_s$ are given by a positive number $P=m^2>0$ and a discrete parameter $s=0, 1/2,1,\dots$ and $p^0 >0$. They correspond to particles with mass $m$ and spin $s$. There is also a class $P^-_s$ with $p^0<0$. The class $0^+_s$ contains representations which correspond to massless particles with discrete helicity and there is also a class $0^-_s$. The representations which are given by  $P=0,\, p^0>0$ and a positive integer $\Xi$ constitute the classes $0^+(\Xi)$, $0^+(\Xi')$ and correspond to elementary systems of mass zero and continuous spin. They are single respectively two valued. The classes $0^-(\Xi)$, $0^-(\Xi)$ are characterized analogously. The remaining classes are given by the cases $p=0$ and $P<0$. 
\end{theorem}
Physical realizations are only known of the classes $P^+_s$ and $0^+_s$. They are realized in Hilbert spaces of $\Lint^2$--integrable functions $\varphi(p,\sigma)$ on the pseudo--riemannian space forms 
\begin{displaymath}
 \Gamma^+_P=\mklm{p \in \R^4: p^kp_k=P,\, p^0 >0}.
\end{displaymath}
The variable $\sigma$ is discrete and can assume the values $-s,\dots.+s$. To each Lorentz transformation $y^k=\Lambda^k{}_lx^l+a^k$ corresponds a unitary operator $U(L)=T(a)d(\Lambda)$ whose action is given by
\begin{equation}
\label{eq:0.0}
 U(L)\varphi(p,\sigma)=e^{i\mklm{p,a}}Q(p,\Lambda) \varphi(\Lambda^{-1}p,\sigma),
\end{equation}
where $Q(p,\Lambda)$ is a unitary operator which depends on $p$ but acts only on the variable $\sigma$. By the continuity of the representations there is for each one--parameter group of unitary operators $U(t)$ an uniquely determined  self adjoint operator $H$ such that $U(t)=\exp(-itH)$.  

In the following we will consider analytic elements of the given representations and use them to characterize real subspaces of localized wave functions. As a corollary of Nelson's analytic vector theorem we have the following proposition \cite{reed-simon}.  
\begin{proposition}
\label{prop:0.0}
A closed  symmetric operator $H$ with domain $D(H)$ acting on a Hilbert space $\H$ is self adjoint if and only if there is a dense set of analytic elements in $D(H)$. The vector--valued functions
  \begin{displaymath}
    U(\tau) \psi:=e^{-i\tau H} \psi= \sum\limits_{n=0}^\infty
    \frac{(-i\tau)^n}{n!} H^n \psi \in \H
  \end{displaymath}
are then analytic in $\tau \in \C$ for each analytic element $\psi$. 
\end{proposition}
The domain of the closed operator $U(\tau)$ depends only on $\Im \tau$. Set  $\tau=\lambda+i\varrho$
and let  $D_U(\varrho)$ be the subset of $\H$ such that $(U(\tau),D_U(\varrho))$ is a closed and normal operator. Then the following statement holds \cite{bisognano-wichmann}.
\begin{proposition}
\label{prop:0.1}
If $\varphi \in D_U(\varrho)$ then the vector--valued function
\begin{displaymath}
  U(\tau)\varphi 
\end{displaymath}
is strongly continuous for $0\leq \Im \tau/\varrho \leq 1$ and analytic for $0<\Im
\tau/\varrho <1$.
\end{proposition}

\section{Boundary values of analytic functions and Fourier--Laplace trans\-form}

We consider boundary values of analytic functions in the sense of distribution theory and specify the necessary conditions for the existence of such limits. We also state the theorems of Paley--Wiener--Schwartz and H\"{o}rmander which will be of relevance in the ensuing sections.    

The following is a generalization of a theorem proved by Epstein \cite{epstein} for boundary values of analytic functions in $\Schw'(\R)$.

\begin{theorem} 
\label{thm:1.0}
Let $\Gamma$ be an open convex cone in $\R^n$ and $\T=\R^n+i\Gamma$. If $f(\zeta$) is analytic in $\T$ and converges for $\Im \zeta \to 0$ to a tempered distribution, that is $\lim _{\Im \zeta \to 0}f(.+i \Im \zeta)$ exists in $\Schw'(\R^n)$, then for each compact set $M$ in $\Gamma$ there is an estimate 
\begin{equation}
  \label{eq:1.1}
  \modulus{f(\zeta)} \leq C(1+\modulus{\zeta})^N, \qquad \Im \zeta \in M,
\end{equation}
where $\modulus{\zeta}:=\max_j\mklm{\zeta_j}$.
\end{theorem}
\begin{proof}
We choose in $\Gamma$ an open convex cone $\Delta$ and $n$ affine independent vectors in $\overline \Delta$ such that in this basis the components of each vector in $\R^n+i\Delta$ have strictly positive imaginary parts. Let $f(.+i\eta)$  converge in $\Schw'(\R^n)$ to the tempered distribution $u$. This means that we can choose $n$ positive numbers $0<\gamma_j<\infty,\, j=1,\dots, n,$ such that for $\eta \in \Delta_\gamma:=\mklm {\eta \in \Delta: 0 < \eta_j < \gamma_j}$ and any $\varphi \in \Schw(\R^n)$ the relation
\begin{displaymath}
  \lim \limits_{\eta \to 0} \int f(\xi+i\eta) \varphi(\xi)
  d\xi=u(\varphi)
\end{displaymath}
holds. We set $f_0=u$ and write this in what follows also as 
\begin{displaymath}
  f_\eta(\varphi)=\eklm{f(.+i\eta),\varphi} \longrightarrow \eklm{f_0,\varphi}=f_0(\varphi).
\end{displaymath}
Since $\Delta$ can be chosen in such a way that $\overline \Delta \setminus\mklm{0} \subset \Gamma$, we obtain a continuous map $\eta \mapsto f_\eta$ from the compactum $\overline \Delta_\gamma$ to $\Schw'(\R^n)$. The image of this map is also compact and hence bounded. One can then infer the existence of a seminorm $\norm{\cdot}_k$ in $\Schw(\R^n)$ such that all $f_\eta$ with $\eta \in \overline {\Delta}_\gamma$ are uniformly continuous with respect to this norm, i.e. there are is a constant $C$ and a positive number $k$ such that for each $\varphi \in \Schw(\R^n)$ and $\eta \in \overline{\Delta}_\gamma$ the inequality
\begin{equation}
  \label{eq:1.2}
  \modulus{\eklm{f_\eta,\varphi}}\leq
  C\sum\limits_{\modulus{\alpha},\,\modulus{\beta} \leq
  k}\sup_x\modulus {x^\alpha \partial ^\beta \varphi(x)}=C\norm{\varphi}_k
\end{equation}
is fulfilled. In this way the distributions $f_\eta$ can be extended to linear
functionals on the Banach space whose topology is induced by the seminorm
$\norm{\cdot}_k$. So the $f_\eta$ can be considered as elements of
${\Schw'}^k(\R^n)$. 

Since $f$ is an analytic function in $\T_\gamma=\R^n+i\Delta_\gamma$, $f(\zeta)/\prod_{j=1}^n(\zeta_j+i)^{-k}$ is also analytic in $\T_\gamma=\R^n+i\Delta_\gamma$ and by Cauchy's integral formula we obtain for $f(\zeta)$ the representation
\begin{displaymath}
  \frac{\prod_{j=1}^n (\zeta_j+i)^k}{(2 \pi i)^n}\int\limits_
  {\kappa_1} \cdots \int\limits_{\kappa_n} \frac {f(z_1,\dots,z_n)}{\prod
  _{j=1}^n (z_j-\zeta_j)(z_j+i)^k} dz_1\cdots dz_n,
\end{displaymath}
where each $\kappa_j$ is a closed path in the strip $\Omega_j=\mklm{z_j \in \C: 0<\Im z_j<\gamma_j}$ around $\zeta_j$ respectively. This representation is independent of the chosen paths and we may thus take them as borders of the rectangles $[-x^1_j,x^1_j]\times[y^1_j,y^2_j], \, x^1_j>0,\,
0<y^1_j<y^2_j<\gamma_i$. 
If we now let the rectangles approach the strips in which they are contained, the integrals along the borders $[y^1_j,y^2_j]$ disappear and we obtain for $f(\zeta)$ the expression
\begin{displaymath}
  \frac{\prod_{j=1}^n (\zeta_j+i)^k}{(2 \pi i)^n}\sum
  \limits_\theta \pm \int \limits_{-\infty}^\infty \frac
  {f(x+i\theta)}{\prod_{j=1}^n
  (x_j+i\theta_j-\zeta_j)(x_j+i\theta_j+i)^k} dx_1 \cdots dx_n
\end{displaymath}
where $\theta=(\theta_1,\dots,\theta_n)$ and the $\theta_j$ are equal to $y^1_j$  or equal to $\gamma_j$ respectively and the sign in the sum over $\theta$ depends on the orientation of the corresponding borders. We put 
\begin{displaymath}
  \Phi_{\zeta,\theta}(x)=\prod\limits_{j=1}^n
  (x_j+i\theta_j-\zeta_j)^{-1}(x_j+i\theta_j+i)^{-k}.
\end{displaymath}
The functions $\Phi_{\zeta,\theta}$ as well as their derivatives up to order $\leq k$ decrease at infinity faster than any polynomial 
$x^\alpha,\, \modulus {\alpha}\leq k$. On the other hand $\theta\mapsto \Phi_{\zeta,\theta}$ is a continuous map for $\theta_j\not=\eta_j=\Im \zeta_j$ into the Banach space ${\Schw}^k(\R^n)$ and  $\Phi_{\zeta,\theta}$ remains in a compact set for $\theta \to 0$ and fixed $\eta_j$. For $\theta \to 0$ we have therefore
\begin{displaymath}
  \eklm {f_\theta,\Phi_{\zeta,\theta}} \rightarrow \eklm{f_{0},\Phi_{\zeta,0}}.
\end{displaymath}
We obtain for $f$ in $\T_\gamma$ the expression
\begin{displaymath}
  f(\zeta)=\frac{\prod_{j+1}^n (\zeta_j+i)^k}{(2 \pi i)^n}\sum
  \limits_\theta \pm \eklm {f_\theta,\Phi_{\zeta,\theta}},
\end{displaymath}
where $\theta=(\theta_1,\dots,\theta_n)$ is equal to $\gamma$ or zero.

As a consequence of the continuity condition \eqref{eq:1.2} it follows for $f(\zeta)$ the estimate
\begin{align}
  \label{eq:1.3}
    \begin{split} 
      \modulus{f(\zeta)} &\leq C\prod \limits_{j=1}^n \modulus{\zeta_j+i}^k
      \sum\limits_\theta \norm {\Phi_{\zeta,\theta}}_k\leq \\
      &\leq C\prod\limits_{j=1}^n(1+\modulus{\zeta_j})^k \sum\limits_\theta
      \sum \limits_{\modulus{\alpha}, \, \modulus{\beta} \leq k} \sup
      _x \modulus{x^\alpha \partial ^\beta \Phi_{\zeta,\theta}(x)},
    \end{split}
\end{align}
where we have made use of
\begin{displaymath}
  \modulus{(\zeta_j+i)^k} \leq \sum \limits_{m=0}^k\modulus{\zeta_j^m
  i^{k-m}}=\sum \limits_{m=0}^k \modulus{\zeta^m_j}=(1+\modulus {\zeta_j})^k.
\end{displaymath}
Differentiation yields for the sum over 
$\modulus{\alpha},\modulus{\beta}\leq k$
\begin{gather*}
  \sum\limits_{\modulus{\alpha},\modulus{\beta}\leq k \atop
  \delta_j+\epsilon_j=\beta_j} \sup _x \modulus{x^\alpha
  \frac{\delta!k!}{\epsilon!} (x_1+i\theta_1-\zeta_1)^{-i-\delta_1}\cdots
  (x_n+i\theta_n-\zeta_n)^{-1-\delta_n} \times \right.\\
  \left.  \times (x_1+i\theta_1+i)^{-k-\epsilon_1}\cdots (x_n+i\theta_n+i)^{-k-\epsilon_n}}
\end{gather*}
so that the estimate \eqref{eq:1.3} now reads
\begin{equation}
  \label{eq:1.4}
  \begin{split}
    \modulus{f(\zeta)} &\leq C\prod \limits_{j=1}^n (1+\modulus{\zeta_j})^k
      \sum\limits_\theta \sum \limits_{\modulus{\delta} \leq k}
      \modulus {\Im  \zeta_1-\theta_1}^{-\delta_1-1}\cdots
      \modulus{\Im \zeta_n-\theta_n}^{-\delta_n-1}.
  \end{split}
\end{equation}
The $\theta_j$ are equal zero or equal $\gamma_j$ and we have $0<\Im \zeta_j<\gamma_j$. Since $\Delta$ and the numbers $\gamma_j$ are arbitrary, we obtain for each compact  $M$ in $\Gamma$ an estimate of the form
\begin{displaymath}
  \modulus{f(\zeta)} \leq C(1+\modulus{\zeta})^N, \qquad \Im \zeta \in M,
\end{displaymath}
for a positive integer $N$ and a constant $C$.
\end{proof}

It can be shown that the given conditions, especially relations \eqref{eq:1.3} and  \eqref{eq:1.4}, are also sufficient.

The Fourier transform is an isomorphism of $\Schw$, so that the Fourier transform of a tempered Distribution $u \in \Schw'(\R^n)$ can be defined as
\begin{displaymath}
\hat u(\varphi)=u(\hat \varphi), \qquad \varphi \in \Schw.
\end{displaymath}
For distributions with compact support the Fourier transform is given by the entire analytic function
\begin{displaymath}
  \hat u(\zeta)=u_x(e^{-i\eklm{x,\zeta}}).\qquad \zeta \in \C^n.
\end{displaymath}
It is called the Fourier--Laplace transform of $u$.

For each closed convex set $E$ we define now the convex, positively homogeneous function
 \begin{equation}
  \label{eq:1.5}
  H_E(\xi):=\sup _{x \in E} \eklm {x,\xi}, \quad \xi \in \R^n.
\end{equation}
with values in $(-\infty,\infty]$.
It characterizes the set $E$ completely, since $E$ is given as the set of all $x \in \R^n$ for which $\eklm{x.\xi} \leq H_E(\xi), \, \xi \in \R^n$. Conversely, if $H$ is a function with the mentioned properties, there exists exactly one closed convex set $E$ such that $H=H_E$ and  $E=\mklm{x: \eklm{x,\xi} \leq H(\xi),\, \xi \in \R^n}$.
If $E$ is compact then $H_E(\xi)<\infty$ for each $\xi$.

We state now the theorems of Paley--Wiener--Schwartz \cite{schwartz} and H\"{o}r\-mander \cite{hoermander}.

\begin{theorem}
\label{thm:1.1}
Let $K$ be a compact convex set in $\R^n$ with support function
$H_K$. If  $u$ is a distribution of order $N$ with support contained in $K$, then 
\begin{equation}
  \label{eq:1.6}
  \modulus {\hat u(\zeta)} \leq C(1+\modulus{\zeta})^N e^{H_K(Im
  \zeta)},\qquad \zeta \in \C^n.
\end{equation}
Conversely, every entire analytic function in $\C^n$ which satisfies the relation \eqref{eq:1.6} for some $N$ is the Fourier--Laplace transform of a distribution with support in $K$.
\end{theorem}

It turns out that it is possible to define the Fourier--Laplace transform at least on certain subspaces of $\C^n$ for more general distributions. So for $\zeta \in \C^n$ and fixed $\eta=\Im \zeta$
\begin{displaymath}
  \hat u(\zeta)= \eklm {u,e^{-i\eklm {.\,,\zeta}}}
\end{displaymath}
could be defined as a distribution in $\xi=\Re \zeta$ if $e^{\eklm{.,\eta}}u \in \Schw'$. We set
\begin{equation}
  \label{eq:1.7}
  \Gamma_u=\mklm{\eta \in \R^n: e^{\eklm{.\,,\eta}} u \in \Schw'}.
\end{equation}
Then the following theorem holds.

\begin{theorem}
\label{thm:1.2}
If  $u \in \D'(\R^n)$, \eqref{eq:1.7} defines a convex set $\Gamma_u$. If its interior $\Gamma_u^\circ$ is not empty, there exists a function $\hat u$
analytic in $\R^n+i \Gamma^\circ_u$ such that the Fourier transform of $e^{\eklm{.,\eta}} u$ is given by $\hat u(.+i \eta)$ for all $\eta \in \Gamma ^\circ_u$. For each compact set $M \subset \Gamma ^\circ_u$ there is an estimate
\begin{equation}
\label{eq:1.8}
  \modulus {\hat u(\zeta)} \leq C (1+ \modulus {\zeta})^N, \qquad \Im
  \zeta \in M.
\end{equation}
Conversely, if $\Gamma$ is an open convex set in $\R^n$ and $U$ an analytic function in $\R^n+i\Gamma$ which fulfills an estimate of the form  \eqref{eq:1.8} for every compact set $M$ in $\Gamma$, then there is a distribution $u$ such that $e^{\eklm{.,\eta}} u \in
\Schw'$  with Fourier transform $U(\cdot+i \eta)$ for all $\eta \in \Gamma$.

If in addition $\supp u \subset K$, then
\begin{equation}
\label{eq:1.9}
  \modulus {\hat u(\zeta)} \leq C (1+\modulus {\zeta})^N e ^{H_K(\Im
  \zeta -\eta)},
\end{equation}
if $\eta \in M$, $H_K(\Im \zeta - \eta) < \infty$, where $M$ is a compact set in $\Gamma^\circ_u$. If conversely there is an $\eta$ for which \eqref{eq:1.9} holds, then $\supp u \subset K$ if $K$ is closed and convex.
\end{theorem}

\section{Localization for the massive scalar field}

Algebraic quantum field theory is concerned with von Neumann algebras $\M(\O)$ of observables localized in space time domains $\O$ together with states $\omega$ on these algebras satisfying certain selection criteria. Due to the Reeh--Schlieder property of the vacuum one may associate with certain regions $\O$ and states $\omega$ the operators $\delta$ and $j$ of the modular theory of Tomita and Takesaki. Here $\delta$ is a positive operator which generates a one parameter group of automorphisms $\ad \delta^{it}$ of $\M(\O)$ and $j$ is an antiunitary operator that defines the conjugation $\ad j$ which maps $\M(\O)$ onto its commutant in the Hilbert space associated with $\omega$ by the Gelfand--Neumann-Segal construction.

Important features of the theory are contained in these operators but explicit realizations of them are only known for certain regions, $\omega$ being the vacuum state. So in the case where $\O$ is a spacelike wedge and the local algebras are generated by Wightman fields that transform covariantly under a  finite dimensional representation of the Lorentz group the modular group is the group of velocity transforms that leave the wedge invariant, and the conjugation is the $PCT$ operation combined with a rotation.

With the knowledge of these modular objects for wedge like regions we associate, following Schroer \cite{schroer1}, to certain closed and convex sets in Minkowski space which arise out of the intersection of wedges real subspaces of wave functions in Wigner representation space. These wave functions can then be viewed as localized in the corresponding regions. 

We characterize the elements of these subspaces as boundary values of analytic functions on three--dimensional complex submanifolds of complex Minkowski space which satisfy certain boundary conditions and show that the latter can be analytically continued to open regions in Minkowski space. They converge in the sense of distribution theory to square--integrable functions and we show that they satisfy the sufficient conditions of the theorems of Paley--Wiener--Schwartz and H\"{o}rmander.

In the following we will restrict ourselves to the massive scalar field. This field corresponds to the representations of class $\P^+_0$ of the inhomogeneous Lorentz group $\Linh$. The wave functions have only one component and the unitary operator  $Q(p,\Lambda)$ in equation \eqref{eq:0.0} is equal to $1$. In this case the wave equation reduces to $p^kp_k=P$. To each Lorentz transformation $y^k=\Lambda^k_lx^l+a^k$ corresponds a unitary operator $U(L)=T(a)d(\Lambda)$ whose action  on any $\varphi \in \Lint^2(\hypP)$ is given by
\begin{equation}
\label{eq:2.0}
  U(L)\varphi(p)=e^{i\mklm{p,a}}\varphi(\Lambda^{-1}p),
\end{equation}
while the $PCT$ transformation is realized by the antiunitary operator
\begin{displaymath}
  \Theta \varphi(p)=\bar \varphi(p).
\end{displaymath}

We consider now in Minkowski space the region $W:=\mklm{x \in \R^4:x^3>\modulus {x^0}}$. $W$ is open and convex as well as invariant under velocity transformations in $x^3$-- direction, under rotations around the $x^3$-- axis and under translations in direction of $x^1$ and $x^2$. All these transformations constitute a subgroup of isometric isomorphisms of $ W$ in $\Linh$. The velocity transforms in $x^3$-- direction 
\begin{displaymath}
  y^k=\Lambda^k{}_l(t)x^l
\end{displaymath}
are given in their active form in the coordinates $x^0, x^1, x^2, x^3$ by the matrices
\begin{displaymath}
  \left ( \begin{array}{cccc} \cosh t & 0 & 0 & -\sinh t \\ 0 & 1 & 0
  & 0 \\ 0 & 0 & 1 & 0 \\ -\sinh t & 0 & 0 & \cosh t \end{array}
  \right ), \qquad t \in \R\,.
\end{displaymath}
By transformation with elements $L=(a,I)$ of $\Linh$ we obtain from $W$ open convex regions $W_L$ to which correspond again certain subgroups of isometric automorphisms. The corresponding velocity transformations
\begin{displaymath}
  y^k=I^k{}_l\Lambda^l{}_m(t)(I^{-1})^m{}_nx^n+(\delta ^k_n-I^k{}_l\Lambda^l{}_m(t)(I^{-1})^m{}_n)a^n
\end{displaymath}
constitute a one--parameter abelian subgroup to which we associate the one--parame\-ter group of unitary operators
\begin{equation}
  \label{eq:2.2}
  U_L(t):=U(({\bf 1}-I\Lambda(t) I^{-1})a,\tilde I \tilde \Lambda(t) \tilde
  I^{-1})=T(a)d(\tilde I) d(\tilde \Lambda(t))d(\tilde I)^{-1}T(a)^{-1}.
\end{equation}
Finally, each $L$ defines the antiunitary involution 
\begin{equation}
  \label{eq:2.3}
  j_L:=T(a)d(\tilde I) d(\tilde \Upsilon) \Theta d(\tilde
  I)^{-1}T(a)^{-1}=T({\bf 1} + I\Upsilon I^{-1})a)d(\tilde I \tilde \Upsilon
  \tilde I^{-1}) \Theta,
\end{equation}
where $\Upsilon$ denotes a rotation by $\pi$ around the $x^3$--axis.

Following Bisognano and Wichmann \cite{bisognano-wichmann} we consider the analytic continuation of the operators $U_L(t)$ and define the closed operators
\begin{displaymath}
  (\delta^{1/2}_{L,+},D_+):=(U_L(i\pi),D_{U_L}(\pi)), \qquad
  (\delta^{1/2}_{L,-},D_-):=(U_L(-i\pi),D_{U_L}(-\pi)).
\end{displaymath}
Together with the involution $j_L$ they define the antilinear closed operators 
\begin{displaymath}
  (s_{L,+},D_+):=(j_L\delta ^{1/2}_{L,+},D_+), \qquad 
  (s_{L,-},D_-):=(j_L\delta ^{1/2}_{L,-},D_-).
\end{displaymath}

They are then the modular operators corresponding to the region $W_L$.
We consider now in Minkowski space a closed polyhedral region $K$ with vertices $a_i$, $i=1, \dots,n$, which arises out of the intersection of the family
\begin{equation}
\label{eq:2.4}
\mklm{\overline{W}_L}_{L \in X_K}
\end{equation}
of closed convex regions $\overline{W}_L$; $X_K$ is some subset of $\Linh$ depending on $K$. This family is supposed to decompose into $n$ subfamilies
\begin{equation}
\label{eq:2.5}
\mklm{\overline
  {W}_{(a_i,I)}}_{I \in X_{K,a_i}}
\end{equation}
where the $X_{K,a_i}$ are nonempty closed convex $6$--dimensional subsets of $\Linh$ associated to each vertex $a_i$. These assumptions correspond to the prescription that the intersection over the family $\mklm {\overline {W}_L}_{L \in X_K}$ is to be understood as the intersection of all the regions $\overline W_L$ which contain $K$.

\begin{definition}
Let $K$ be a closed convex region as above. We associate to $K$ in $\Lint^2(\hypP)$ the real subspaces
 \begin{displaymath}
  \H^R_{K,\pm}:=\mklm{\varphi \in \Lint^2(\hypP):
  s_{L,\pm}\varphi=\varphi,\quad L \in X_K}.
\end{displaymath}
\end{definition}

Let $\varphi \in \H^R_{K,+}$. Then $\varphi \in 
D_{U_L}(\pi)$ for all $L\in X_K$ and we define for each vertex $a_i$ the functions 
\begin{equation}
  \label{eq:2.6}
  u_{+,a_i}(\zeta):=U_L(\tau)\varphi(p), \quad L=(a_i,I),\quad I  \in X_{K,a_i},
  \quad 0\leq\Im \tau \leq \pi,
\end{equation}
where $p \in \hypP$ and $\zeta =I\Lambda(\tau)^{-1}I^{-1}p$.

\begin{lemma}
For every vertex $a_i$, the set
\begin{equation}
\label{eq:2.7}
  M_{K,a_i}^+:=\mklm{\zeta \in \chypP: \zeta =I\Lambda(\tau)^{-1}I^{-1}p,
  \, p \in \hypP,\, 0 < \Im \tau <\pi,\, I \in X^\circ_{K,a_i}}
\end{equation}
is a complex 3--dimensional submanifold in $\C^4$ and the function $u_{+,a_i}$ is holomorphic on $M^+_{K,a_i}$ and hence uniquely determined. For $0 \leq \Im \tau \leq \pi$, $u_{+,a_i}$ is continuous.
\end{lemma}
\begin{proof}
It can be shown that the sets $M^+_{K,a_i}$ are given by
\begin{displaymath}
  M_{K,a_i}^+=\chypP \cap \T_{K,a_i}^+,
\end{displaymath}
where $\T_{K,a_i}^+=\R^4+i\Gamma_{K,a_i}^+$ and
\begin{gather*}
  \Gamma_{K,a_i}^+:=\mklm{\eta \in \R^4: \eta=I\eta',\quad I\in X^\circ_{K,a_i},\quad
  \eta' \in \Gamma^+ },\\
  \Gamma^+:=\mklm{\eta \in \R^4:\eta^-<0,\quad \eta^+>0,\quad \eta^1=\eta^2=0}.
\end{gather*}
Since $\chypP$ is given as the set of zeros of the analytic function $\mu_(\zeta)=\mklm{\zeta,\zeta}-P$ it is a complex submanifold of $\C^4$ and since $\T^+_{K,a_i}$ is open in $\C^4$, 
so is $M_{K,a_i}^+$ as well. $\LC$ acts transitively on $\chypP$ and the isotropy group of a point $\zeta \in \chypP$ is given by $\SU(3)$, so that
\begin{displaymath}
  \chypP \simeq \LC \Big / \SU(3).
\end{displaymath}
$\chypP$ is thus the orbit of a point $\zeta$ under complex velocity transforms and we can write each $\zeta$ in $M_{K,a_i}^+=\chypP \cap \T_{K,a_i}^+$ as
\begin{displaymath}
  \zeta=I\Lambda^{-1}(\tau_I)I^{-1}p=\Lambda^{-1}_1(\tau_1)\Lambda^{-1}_2(\tau_2)\Lambda^{-1}_3(\tau_3)p_0, 
\end{displaymath}
where  the $\Lambda_i(\tau_i)$ are complex velocity transforms in $x_i$--direction with $0 \leq \Im \tau_i \leq \pi,\, i=1,2,3,$ and
$p_0$ is a fixed point in $\hypP$. We can thus interpret the functions 
$u_{+,a_i}(\zeta)$ as functions of $\tau=(\tau_1,\tau_2,\tau_3)$ and write
$u_{+,a_i}(\zeta)=u_{+,a_i}(\tau_1,\tau_2,\tau_3)$. Since every complex velocity transform can be obtained from any other by adjoining the latter with a homogeneous Lorentz transform, the analycity of the $u_{+,a_i}$ follows by 
Proposition \ref{prop:0.0}.
\end{proof}
The $u_{+,a_i}$ satisfy the boundary conditions
\begin{displaymath}
  u_{+,a_i}(-\xi)=U_L(i\pi)\varphi(p)=j_L\varphi(p)=e^{i\mklm{p,a_i}}e^{i\mklm{\xi,a_i}}\bar
  u_{+,a_i}(\xi) 
\end{displaymath}
for all $\xi=I\Upsilon^{-1}I^{-1}p,\, p \in
\hypP,\,I \in X_{K,a_i}$. Note that $\Lambda^{-1}(i\pi)=-\Upsilon$. Analogous considerations hold in the case that $\varphi\in \H^R_{K,-}$. Hence we have the following proposition.
\begin{proposition}
\label{prop:2.0}
  Let an irreducible representation of the inhomogeneous Lor\-entz group $\Linhgen$ of class $P^+_0$ be given in the Hilbert space $\Lint^2(\hypP)$ of $\Lint^2$--integrable functions on $\hypP$. The subspaces  $\H^R_{K,\pm}$ associated to the polyhedral region $K$ with vertices $a_i$, $i=1,\dots,n$, 
are then given by the set of all functions $\varphi(p) \in \Lint^2(\hypP)$ which are boundary values of analytic functions $u_{+,a_i}$ resp. $u_{-,a_i}$, in $M_{K,a_i}^+$ resp. $M_{K,a_i}^-$ for some $i$ satisfying the boundary conditions  
  \begin{equation}
    \label{eq:2.8}
    u_{\pm,a_i}(-\xi)=e^{i\mklm{p,a_i}}e^{i\mklm{\xi,a_i}}\bar u_{\pm,a_i}(\xi)
  \end{equation}
for all $\xi=I\Upsilon I^{-1}p,\, p \in \hypP$, and $(a_i,I)  \in
\mklm{a_i} \times X_{K,a_i} \subset X_K$ respectively. 
\end{proposition}

By the second theorem of Cartan every holomorphic function on a closed analytic submanifold of a Steinian manifold $X$ is the restriction of a holomorphic function defined on $X$ \cite{gunning-rossi}. Therefore every$u_{+,a_i}$ can be extended analytically to $\T^+_{K,a_i}$ respectively. Since the holomorphic hull of $\T^+_{K,a_i}$ is given by its convex hull, we have the following lemma.

\begin{lemma}
$u_{+,a_i}$ can be analytically extended to $\convh \T^+_{K,a_i}$.
\end{lemma}
 
We consider now the operator $R(\Lambda)\varphi(p)=\varphi(\Lambda^{-1}p)$ and set 
\begin{displaymath}
R_I(t):=R(I \Lambda^{-1}(t) I^{-1}).
\end{displaymath}
By definition we have
\begin{displaymath}
  U_L(t)\varphi(p)=e^{i\mklm{p,a}}e^{-i\mklm{I\Lambda^{-1}(t)I^{-1}p,a}}R_I(t)\varphi(p).
\end{displaymath}
If $\varphi \in D_R(\pi)$ we define the function
\begin{equation}
  \label{eq:2.9}
  r_+(\zeta):=R_I(\tau)\varphi(p), \quad (a,I) \in X_K,\quad 0\leq\Im \tau\leq \pi.
\end{equation} 
In a similar way as before for the functions $u_{+,a_i}$ we have the following lemma.
\begin{lemma}
$r_+$ is holomorphic  on the complex manifold
\begin{displaymath}
  M_{K,+}:=\mklm{\zeta \in \chypP: \zeta =I\Lambda(\tau)^{-1}I^{-1}p,
  \, p \in \hypP,\, 0 < \Im \tau <\pi,\, L=(a,I) \in X_{K}^\circ}
\end{displaymath}
and can analytically be extended to $\convh \T_{K}^+$, where  
\begin{displaymath}
\Gamma_K^+:=\mklm{\eta \in \R^4: \eta=I\eta',\, (a,I)\in X^\circ_{K},\,
  \eta' \in \Gamma^+ }. 
\end{displaymath}
\end{lemma}

As before we can represent each element $\zeta$ of $ M_{K,+}=\chypP \cap \T^+_K$ as
\begin{displaymath}
  \zeta=I\Lambda^{-1}(\tau_I)I^{-1}p=\Lambda^{-1}_1(\tau_1)\Lambda^{-1}_2(\tau_2)\Lambda^{-1}_3(\tau_3)p_0, 
\end{displaymath}
where  $\Lambda_i(\tau_i)$ are complex velocity transformations in $x_i$ direction with $0 \leq \Im \tau_i \leq \pi,\, i=1,2,3,$ and
$p_0$ a fixed point in $\hypP$. We can therefore interpret $r_+(\zeta)$ as a function of $\tau=(\tau_1,\tau_2,\tau_3)$ and write
\begin{gather*}
  r_+(\zeta)=r_+(\tau_1,\tau_2,\tau_3)=\\
  =R(\Lambda_3(\tau_3))R(\Lambda_2(\tau_2))R(\Lambda_1(\tau_1))\varphi(p_0)
  =R_{\Lambda_3\Lambda_2\Lambda_1}(\tau).
\end{gather*}

\begin{proposition}
$r_+$ satisfies in $\convh \T^+_K$ an estimate of the form
\begin{equation}
  \label{eq:2.10}
  \modulus{r_+(\zeta)} \leq C(1+\modulus{\zeta})^N,\qquad \Im \zeta
  \in H,
\end{equation}
where $H$ is a compact set in  $\convh \Gamma^+_K$.
\end{proposition}
\begin{proof}
The element $R_I(\tau_I)\varphi(.)$in $\Lint^2(\hypP)$  is a tempered distribution 
\begin{displaymath}
R_{\Lambda_3\Lambda_2\Lambda_1}(\tau)\varphi(.)=
R_{\Lambda_3\Lambda_2\Lambda_1}(.+i\varrho)\varphi(p_0)
\end{displaymath}
 which depends on the parameters $\varrho=(\varrho_1,\varrho_2,\varphi_3)$; for $\Phi \in \Schw(\hypP)$ we have thus 
\begin{gather*}
  \eklm {R_{\Lambda_3\Lambda_2\Lambda_1}(i\varrho)\varphi,\Phi}=
  \\=\int^\infty_{-\infty}R_{\Lambda_3\Lambda_2\Lambda_1}(\lambda'+i\varrho)\varphi(p_0)
  \Phi(\Lambda_1^{-1}(\lambda_1')\Lambda_2^{-1}(\lambda_2')\Lambda_3^{-1}(\lambda_3')p_0)\,
  d\lambda'=\\
  =\int_\hypP R_{\Lambda_3\Lambda_2\Lambda_1}(i\varrho)\varphi(p)\Phi(p) \, dM
\end{gather*}
where $\modulus {\partial (p^1,p^2,p^3)/\partial(\lambda_1,\lambda_2,\lambda_3)}=p^0$. 
Now $R_I(\tau_I)\varphi=R_{\Lambda_3\Lambda_2\Lambda_1}(\tau)\varphi$
is strongly continuous for  $0 \leq \Im \tau_I \leq \pi$ as well as for corresponding values of $\tau_i,\, i=1,2,3$; in particular we obtain
\begin{displaymath}
  \norm {R_{\Lambda_3\Lambda_2\Lambda_1}(\tau)\varphi-\varphi} \to 0
  \quad \text{for} \quad \tau \to 0;
\end{displaymath}
since strong convergence implies weak convergence we have for all elements $\psi$ of $\Lint^2(\hypP)$
\begin{displaymath}
  (R_{\Lambda_3\Lambda_2\Lambda_1}(\tau)\varphi,\psi) \to
  (\varphi,\psi) \quad \text{for} \quad \tau \to 0.
\end{displaymath}
Hence it follows that
\begin{displaymath}
  \lim_{\varrho \to 0} \int_\hypP R_{\Lambda_3\Lambda_2\Lambda_1}(i\varrho)\varphi(p)\Phi(p)\,dM=\varphi(\Phi),
\end{displaymath}
and by Theorem \ref{thm:1.0} we obtain for  $r_+(\zeta)=R_{\Lambda_3\Lambda_2\Lambda_1}(\tau)\varphi(p_0)$ in $M_{K,+}$ an estimate of the form
\begin{equation}
  \label{eq:2.10.a}
  \modulus{r_+(\zeta)} \leq C(1+\modulus{\zeta})^N,\qquad \Im \zeta
  \in H,
\end{equation}
where $H$ is a compact set in $\Im M_{K,+}$. Continuing
$r_+\equiv r_+(\zeta^i),\, \zeta^i\equiv \tau_i,\, \zeta^4=\mu(\zeta)-P$, analytically to $\convh
\T^+_K$ in such a way that $r_+$ also has a bound of the given form with respect to $\zeta^4$ we obtain for $r_+$ in all $\convh \T^+_K$ an estimate of the form \eqref{eq:2.10.a} where  $H$
is now a compact set in $\convh \Gamma^+_K$.
\end{proof}
For all functions $u_{+,a_i}$ the relation
\begin{equation}
\label{eq:2.11}
  u_{+,a_i}(\zeta)=e^{i\mklm{p,a_i}}e^{-i\mklm {\zeta,a_i}}r_+(\zeta),\qquad
  \zeta \in \convh \T_{K,a_i}^+
\end{equation}
holds and we obtain with \eqref{eq:2.10} the estimates
\begin{align*}
  \modulus{u_{+,a_i}(\zeta)}&=e^{\mklm{\Im \zeta,a_i}}\modulus {r_+(\zeta)}\leq\\
  &\leq C(1+\modulus{\zeta})^N e^{H_K^{(1,3)}(\Im \zeta)}, \quad \Im
  \zeta \in H,
\end{align*}
for each compact set $H$ in $\convh \Gamma_{K,a_i}^+$ and all 
$a_i$, $i=1,\dots,n$. We reformulate this and summarize the above results in the following theorem.
\begin{theorem}
\label{thm:2.0}
 Let $K$ be a polyhedral region in $\R^4$ with vertices $a_i$,
 $i=1,\dots,n$, and $\varphi \in \H^R_{K,+}$ resp. $\H^R_{K,-}$ the boundary value of the functions $u_{+,a_i}$ resp. $u_{-,a_i}$ analytic in $M_{K,a_i}^+=\chypP \cap \T_{K,a_i^+}$ resp. $M_{K,a_i}^-=\chypP \cap \T_{K,a_i}^-$. Then $u_{+,a_i}$ and
 $u_{-,a_i}$ can be analytically extended to $\convh \T_{K,a_i}^+$ resp. $\convh \T_{K,a_i}^-$ and satisfy for each compact set $H$ in $\convh \Gamma_{K,a_i}^+$
 resp. $\convh \Gamma_{K,a_i}^-$ an estimate
 \begin{equation*}
   \label{eq:2.12}
   \modulus{u_{\pm,a_i}(\zeta)} \leq
   C(1+\modulus{\zeta})^Ne^{H_{K}^{(1,3)}(\Im \zeta-\eta)},\qquad
   \zeta \in \convh \T_{K,a_i}^\pm,
 \end{equation*}
if $\eta \in H$ and $H_{K}^{(1,3)}(\Im \zeta -\eta) < \infty $.
\end{theorem}
Thus the functions $u_{+,a_i}$ and $u_{+,a_i}$ satisfy the sufficient conditions of Theorem \ref{thm:1.2}. If in addition $K$ is compact we have the following theorem.
\begin{theorem}
\label{thm:2.1}
  Let $K$ a compact polyhedral region in $\R^4$ and $\varphi \in \H^R_{K,+}$ resp. $\H^R_{K,-}$.
  Then $\varphi$ is the boundary value of an entire analytic function $u(\zeta)$ which satisfies an estimate of the form
   \begin{displaymath}
     \modulus{u(\zeta)} \leq C(1+\modulus{\zeta})^N e^{H_K^{(1,3)}(\Im
     \zeta)}, \qquad \zeta \in \C^4. 
   \end{displaymath}
\end{theorem}

According to Theorem \ref{thm:1.1} the elements of $\H^R_{K,+}$ resp. $\H^R_{K,-}$ are then boundary values of Fourier--Laplace transforms of distributions with support in $K$.


\providecommand{\bysame}{\leavevmode\hbox to3em{\hrulefill}\thinspace}



\begin{thebibliography}{Tom67}

\bibitem[BGLa]{brunetti-guido-longo2}
R.~Brunetti, D.~Guido, and R.~Longo, \emph{First quantization via the
  {{\rm{BW}}}--property}, in progress.

\bibitem[BGLb]{brunetti-guido-longo1}
R.~Brunetti, D.~Guido, and R.~Longo, \emph{On the intinsic construction of free
  theories via {Tomita--Takesaki}--theory}, unpublished manuscript.

\bibitem[BW48]{bargmann-wigner}
V.~Bargmann and E.~P. Wigner, \emph{Group theoretical discussion of
  relativistic wave equations}, Proc. Nat. Acad. Sc. USA \textbf{34} (1948),
  211--223.

\bibitem[BW75]{bisognano-wichmann}
J.~J. Bisognano and E.H. Wichmann, \emph{On the duality condition for a
  {Hermitian} scalar field}, Jour. Math. Phys. \textbf{16} (1975), 985--1007.

\bibitem[Eps66]{epstein}
H.~Epstein, \emph{Some analytic properties of scattering amplitudes in quantum
  field theory}, Axiomatic Field Theory, Brandeis University Summer Institute
  in Theoretical Physics, 1965, vol.~1, Gordon and Breach, Science Publishers,
  Inc., New York, 1966, pp.~1--133.

\bibitem[GR65]{gunning-rossi}
R.C. Gunning and H.~Rossi, \emph{Analytic functions of several complex
  variables}, Prentice--Hall, INC., Englewood Cliffs, New York, 1965.

\bibitem[Haa92]{haag}
R.~Haag, \emph{Local quantum physics, fields, particles, algebras},
  Springer--Verlag, Berlin, Heidelberg, New York, 1992.

\bibitem[H{\"{o}}r83]{hoermander}
L.~H{\"{o}}rmander, \emph{The analysis of linear partial differential
  operators}, vol.~I, Springer--Verlag, Berlin, Heidelberg, New York, 1983.

\bibitem[RS75]{reed-simon}
M.~Reed and B.~Simon, \emph{Methods of modern mathematical physics}, vol.~II,
  Academic Press, INC. San Diego, London, 1975.

\bibitem[Sch51]{schwartz}
L.~Schwartz, \emph{Th\'{e}orie de distributions}, vol.~II, Hermann \&
  $\rm{C^{ie}}$, Paris, 1951.

\bibitem[Sch97]{schroer1}
B.~Schroer, \emph{Wigner representation theory of the {Poincar\'{e}} group,
  localization, statistics and the {S--matrix}}, Nucl. Phys. B \textbf{499}
  (1997), 519--546.

\bibitem[Sch99]{schroer3}
B.~Schroer, \emph{Modular wedge localization and the $d=1+1$ formfactor
  program}, Annals of Physics \textbf{275} (1999), 190--223.

\bibitem[Tak70]{takesaki}
M.~Takesaki, \emph{Tomita's theory of modular {Hilbert}--algebras and its
  application}, Lecture Notes in Mathematics, vol. 128, Springer--Verlag,
  Berlin, Heidelberg, New York, 1970.

\bibitem[Tom67]{tomita}
M.~Tomita, \emph{Quasi--standard von {Neumann} algebras}, mimeographed notes,
  1967.

\bibitem[Wig39]{wigner1}
E.~P. Wigner, \emph{On unitary representations of the inhomogenous {Lorentz}
  group}, Ann. Math. \textbf{40} (1939), no.~1, 149--204.

\end{thebibliography}
\end{document}